\input ./style/arxiv-general.cfg
\documentclass[MSNbibl,number,citesort,seceqn,dvips]{arxbj}
\makeatletter
   \@ifpackageloaded{graphicx}{}{\usepackage{graphicx}}
\makeatother

\aid{0}
\volume{21}
\issue{3}
\pubyear{2015}
\firstpage{1800}
\lastpage{1823}
\doi{10.3150/14-BEJ625} 

\makeatletter

\newcommand{\rrvert}{\vert}
\newcommand{\llvert}{\vert}
\newcommand{\overset}{\stackrel}
\newproclaim{definition}{Definition}[section]
\newproclaim{assumption}{Assumption}[section]
\newtheorem{theorem}{Theorem}[section]
\newtheorem{corollary}{Corollary}[section]
\newtheorem{lemma}{Lemma}[section]
\makeatother

\begin{document}
\begin{frontmatter}

\title{Interplay of insurance and financial risks in~a~discrete-time
model with strongly regular~variation}
\runtitle{Interplay of insurance and financial risks}

\begin{aug}
\author[1]{\inits{J.}\fnms{Jinzhu}~\snm{Li}\corref{}\thanksref{1}\ead[label=e1]{lijinzhu@nankai.edu.cn}} \and
\author[2]{\inits{Q.}\fnms{Qihe}~\snm{Tang}\thanksref{2}\ead[label=e2]{qihe-tang@uiowa.edu}}
\address[1]{School of Mathematical Science and LPMC, Nankai University,
Tianjin 300071, P.R. China.\\ \printead{e1}}
\address[2]{Department of Statistics and Actuarial Science, The
University of Iowa, 241 Schaeffer Hall, Iowa City, IA 52242, USA.
\printead{e2}}
\end{aug}

\received{\smonth{1} \syear{2013}}
\revised{\smonth{11} \syear{2013}}

%
\begin{abstract}
Consider an insurance company exposed to a stochastic economic environment
that contains two kinds of risk. The first kind is the insurance risk caused
by traditional insurance claims, and the second kind is the financial risk
resulting from investments. Its wealth process is described in a standard
discrete-time model in which, during each period, the insurance risk is
quantified as a real-valued random variable $X$ equal to the total
amount of
claims less premiums, and the financial risk as a positive random
variable $%
Y $ equal to the reciprocal of the stochastic accumulation factor. This risk
model builds an efficient platform for investigating the interplay of the
two kinds of risk. We focus on the ruin probability and the tail probability
of the aggregate risk amount. Assuming that every convex combination of the
distributions of $X$ and $Y$ is of strongly regular variation, we derive
some precise asymptotic formulas for these probabilities with both finite
and infinite time horizons, all in the form of linear combinations of the
tail probabilities of $X$ and $Y$. Our treatment is unified in the sense
that no dominating relationship between $X$ and $Y$ is required.
\end{abstract}

%
\begin{keyword}
\kwd{asymptotics}
\kwd{convolution equivalence}
\kwd{financial risk}
\kwd{insurance risk}
\kwd{ruin probabilities}
\kwd{(strongly) regular variation}
\kwd{tail probabilities}
\end{keyword}
\end{frontmatter}

\section{Introduction}

As summarized by Norberg \cite{Nor99}, an insurance
company which invests its
wealth in a financial market is exposed to two kinds of risk. The first
kind, called insurance risk, is the traditional liability risk caused by
insurance claims, and the second kind, called financial risk, is the asset
risk related to risky investments. The interplay of the two risks
unavoidably leads to a complicated stochastic structure for the wealth
process of the insurance company. Paulsen \cite{Pau93} proposed a general
continuous-time risk model in which the cash flow of premiums less
claims is
described as a semimartingale and the log price of the investment portfolio
as another semimartingale. Since then the study of ruin in the presence of
risky investments has experienced a vital development in modern risk theory;
some recent works include Paulsen \cite{Pau08}, Kl\"{u}ppelberg and Kostadinova
\cite{KluKos08}, Heyde and Wang \cite{HeyWan09},
Hult and Lindskog \cite{HulLin11}, Bankovsky \textit{et al.}
\cite{BanKluMal11}, and Hao and Tang \cite{HaoTan12}. During this
research, much attention has
been paid to an important special case of Paulsen's set-up, the so-called
bivariate L\'{e}vy-driven risk model, in which the two semimartingales are
independent L\'{e}vy processes fulfilling certain conditions so that
insurance claims dominate financial uncertainties.

A well-known folklore says that risky investments may impair the insurer's
solvency just as severely as do large claims; see Norberg
\cite{Nor99},
Kalashnikov and Norberg \cite{KalNor02}, Frolova \textit{et al.} \cite{FroKabPer02},
and Pergamenshchikov and
Zeitouny \cite{PerZei06}.

In this paper, we describe the insurance business in a discrete-time risk
model in which the two risks are quantified as concrete random variables.
This discrete-time risk model builds an efficient platform for investigating
the interplay of the two risks. The ruin probabilities of this model have
been investigated by
Nyrhinen \cite{Nyr99,Nyr01},
Tang and Tsitsiashvili \cite{TanTsi03,TanTsi04}, Collamore \cite{Col09}, and
Chen \cite{Che11}, among many others.

Concretely, for each $n\in\mathbb{N}=\{1,2,\ldots\}$, denote by $X_{n}$
the insurer's net loss (the total amount of claims less premiums) within
period $n$ and by $Y_{n}$ the stochastic discount factor (the
reciprocal of
the stochastic accumulation factor) over the same time period. Then the
random variables $X_{1}$, $X_{2}$, \ldots\ and $Y_{1}$, $Y_{2}$,
\ldots\
represent the corresponding insurance risks and financial risks,
respectively. In this framework, we consider the stochastic present values
of aggregate net losses specified as
%
%
\begin{equation}
\label{S} S_{0}=0,\qquad S_{n}=\sum
_{i=1}^{n}X_{i}\prod
_{j=1}^{i}Y_{j},\qquad n\in \mathbb{N},
\end{equation}
and consider their maxima
%
%
\begin{equation}
\label{M} M_{n}=\max_{0\leq k\leq n}S_{k},\qquad
n\in\mathbb{N}.
\end{equation}

If $(X_{1},Y_{1})$, $(X_{2},Y_{2})$, \ldots\ form a sequence of independent
and identically distributed (i.i.d.) random pairs fulfilling $-\infty
\leq
\mathrm{E}\ln Y_{1}<0$ and $\mathrm{E}\ln ( \llvert
X_{1}\rrvert
\vee1 ) <\infty$, then, by Lemma~1.7 of Vervaat
\cite{Ver79}, $S_{n}$
converges almost surely (a.s.) as $n\rightarrow\infty$. In this case,
denote by $S_{\infty}$ the a.s. limit. Clearly, $M_{n}$ is non-decreasing
in $n$ and
\[
0\leq M_{n}\leq\sum_{i=1}^{n} (
X_{i}\vee0 ) \prod_{j=1}^{i}Y_{j}.
\]
Thus, if $-\infty\leq\mathrm{E}\ln Y_{1}<0$ and $\mathrm{E}\ln
 (
X_{1}\vee1 ) <\infty$, then $M_{n}$ also converges a.s. to a limit,
denoted by $M_{\infty}$, as $n\rightarrow\infty$.

We conduct risk analysis of the insurance business through studying the tail
probabilities of $S_{n}$ and $M_{n}$ for $n\in\mathbb{N}\cup\{\infty
\}$.
The study of tail probabilities is of fundamental interest in insurance,
finance, and, in particular, quantitative risk management. Moreover, the
tail probability of $M_{n}$ with $n\in\mathbb{N}\cup\{\infty\}$ is
immediately interpreted as the finite-time or infinite-time ruin probability.

In most places of the paper, we restrict ourselves to the standard framework
in which $X_{1}$, $X_{2}$, \ldots\ form a sequence of i.i.d. random
variables with generic random variable $X$ and common distribution $F=1-
\overline{F}$ on $\mathbb{R}=(-\infty,\infty)$, $Y_{1}$, $Y_{2}$,
\ldots\
form another sequence of i.i.d. random variables with generic random
variable $Y$ and common distribution $G$ on $(0,\infty)$, and the two
sequences are mutually independent.

Under the assumption that the insurance risk $X$ has a regularly-varying
tail dominating that of the financial risk $Y$,
Tang and Tsitsiashvili
\cite{TanTsi03,TanTsi04} obtained some precise
asymptotic formulas for the finite-time
and infinite-time ruin probabilities. The dominating relationship
between $X$
and $Y$ holds true if we consider the classical Black--Scholes market in
which the log price of the investment portfolio is modelled as a Brownian
motion with drift and, hence, $Y$ has a lognormal tail, lighter than every
regularly-varying tail. However, empirical data often reveal that the
lognormal model significantly underestimates the financial risk. It shows
particularly poor performance in reflecting financial catastrophes such as
the recent Great Recession since 2008. This intensifies the need to
investigate the opposite case where the financial risk $Y$ has a
regularly-varying tail dominating that of the insurance risk $X$. In this
case, the stochastic quantities in (\ref{S}) and (\ref{M}) become much
harder to tackle with the difficulty in studying the tail probability
of the
product of many independent regularly-varying random variables.
Tang and Tsitsiashvili \cite{TanTsi03} gave two examples for
this opposite case illustrating
that, as anticipated, the finite-time ruin probability is mainly determined
by the financial risk. Chen and Xie \cite{CheXie05}
also studied the finite-time ruin
probability of this model and they obtained some related results applicable
to the case with the same heavy-tailed insurance and financial risks.

In this paper, under certain technical conditions, we give a unified
treatment in the sense that no dominating relationship between the two risks
is required. That is to say, the obtained formulas hold uniformly for the
cases in which the insurance risk $X$ is more heavy-tailed than, less
heavy-tailed than, and equally heavy-tailed as the financial risk $Y$. In
our main result, under the assumption that every convex combination of $F$
and $G$ is of strongly regular variation (see Definition~\ref{SRV} below),
we derive some precise asymptotic formulas for the tail probabilities
of $%
S_{n}$ and $M_{n}$ for $n\in\mathbb{N}\cup\{\infty\}$. All the obtained
formulas appear to be linear combinations of $\overline{F}$ and~$\overline{G}
$. Hence, if one of $\overline{F}$ and $\overline{G}$ dominates the other,
then this term remains in the formulas but the other term is negligible;
otherwise, both terms should simultaneously present. These formulas are in
line with the folklore quoted before, confirming that whichever one of the
insurance and financial risks with a heavier tail plays a dominating
role in
leading to the insurer's insolvency.

In the rest of this paper, Section~\ref{sec2} displays our results and some related
discussions after introducing the assumptions, Section~\ref{sec3} prepares some
necessary lemmas, and Section~\ref{sec4} proves the results.

\section{Preliminaries and results}\label{sec2}

Throughout this paper, all limit relationships hold for $x\rightarrow
\infty
$ unless otherwise stated. For two positive functions $a(\cdot)$ and $%
b(\cdot)$, we write $a(x)\lesssim b(x)$ or $b(x)\gtrsim a(x)$ if
$\limsup
a(x)/b(x)\leq1$, write $a(x)\sim b(x)$ if both $a(x)\lesssim b(x)$ and
$%
a(x)\gtrsim b(x)$, and write $a(x)\asymp b(x)$ if both $a(x)=\mathrm{O}(b(x))$
and $%
b(x)=\mathrm{O}(a(x))$. For a real number $x$, we write $x_{+}=x\vee0$ and $%
x_{-}=-(x\wedge0)$.

\subsection{Assumptions}
\label{A}

We restrict our discussions within the scope of regular variation. A
distribution $U$ on $\mathbb{R}$ is said to be of regular variation if $
\overline{U}(x)>0$ for all $x$ and the relation
\[
\lim_{x\rightarrow\infty}\frac{\overline{U}(xy)}{\overline{U}(x)}%
=y^{-\alpha},\qquad
y>0,
\]
holds for some $0\leq\alpha<\infty$. In this case, we write $U\in
\mathcal{R}_{-\alpha}$. However, such a condition is too general to enable us to
derive explicit asymptotic formulas for the tail probabilities of the
quantities defined in (\ref{S}) and (\ref{M}). To overcome this difficulty,
our idea is to employ some existing results and techniques related to the
well-developed concept of convolution equivalence.

A distribution $V$ on $[0,\infty)$ is said to be convolution
equivalent if
$\overline{V}(x)>0$ for all $x$ and the relations
%
%
\begin{equation}
\label{L(alpha)} \lim_{x\rightarrow\infty}\frac{\overline{V}(x-y)}{\overline
{V}(x)}=\mathrm{%
e}^{\alpha y},
\qquad y\in\mathbb{R},
\end{equation}
and
%
%
\begin{equation}
\label{S(alpha)} \lim_{x\rightarrow\infty}\frac{\overline{V^{2\ast}}(x)}{\overline
{V}(x)}%
=2c<\infty
\end{equation}
hold for some $\alpha\geq0$, where $V^{2\ast}$ stands for the $2$-fold
convolution of $V$. More generally, a distribution $V$ on $\mathbb{R}$ is
still said to be convolution equivalent if $V(x)\mathbf{1}_{(x\geq
0)}$ is.
In this case, we write $V\in\mathcal{S}(\alpha)$. Relation (\ref{L(alpha)})
itself defines a larger class denoted by $\mathcal{L}(\alpha)$. It is known
that the constant $c$ in relation (\ref{S(alpha)}) is equal to
\[
\hat{V}(\alpha)=\int_{-\infty}^{\infty}\mathrm{e}^{\alpha
x}V(
\mathrm {d}%
x)<\infty;
\]
see Cline \cite{Cli87} and Pakes
\cite{Pak04}. We shall use the notation $\hat
{V}(\cdot
) $ as above for the moment generating function of a distribution $V$
throughout the paper. The class $\mathcal{S}(0)$ coincides with the
well-known subexponential class. Examples and criteria for membership
of the
class $\mathcal{S}(\alpha)$ for $\alpha>0$ can be found in
Embrechts \cite{Emb83} and Cline
\cite{Cli86}. Note that the gamma distribution belongs to the
class $\mathcal{L}(\alpha)$ for some $\alpha>0$ but does not belong
to the
class $\mathcal{S}(\alpha)$. Hence, the inclusion $\mathcal
{S}(\alpha
)\subset\mathcal{L}(\alpha)$ is proper. Recent works in risk theory using
convolution equivalence include Kl\"{u}ppelberg \textit{et al.} \cite{KluKypMal04},
Doney and Kyprianou \cite{DonKyp06}, Tang and Wei
\cite{TanWei10}, Griffin and Maller \cite{GriMal12},
Griffin \textit{et al.} \cite{GriMalvan12}, and Griffin \cite{Gri13}.

For a distribution $U$ on $\mathbb{R}$, define
%
%
\begin{equation}
\label{U--V} V(x)=1-\frac{\overline{U}(\mathrm{e}^{x})}{\overline{U}(0)},\qquad x\in \mathbb{R},
\end{equation}
which is still a proper distribution on $\mathbb{R}$. Actually, if
$\xi
$ is
a real-valued random variable distributed as $U$, then $V$ denotes the
conditional distribution of $\ln\xi$ on $\xi>0$. For every $\alpha
\geq0$, it is clear that $U\in\mathcal{R}_{-\alpha}$ if and only if $V\in
\mathcal{L}(\alpha)$. We now introduce a proper subclass of the class $
\mathcal{R}_{-\alpha}$.

%
\begin{definition}
\label{SRV}A distribution $U$ on $\mathbb{R}$ is said to be of strongly
regular variation if $V$ defined by (\ref{U--V}) belongs to the class $
\mathcal{S}(\alpha)$ for some $\alpha\geq0$. In this case, we write
$U\in
\mathcal{R}_{-\alpha}^{\ast}$.
\end{definition}

Examples and criteria for membership of the class $\mathcal
{R}_{-\alpha
}^{\ast}$ can be given completely in parallel with those in
Embrechts \cite{Emb83} and Cline
\cite{Cli86}. This distribution class turns out to be
crucial for
our purpose. Clearly, if $\xi$ follows $U\in\mathcal{R}_{-\alpha
}^{\ast}$
for some $\alpha\geq0$, then
\[
\mathrm{E}\xi_{+}^{\alpha}=\overline{U}(0)\mathrm{E} \bigl(
 \mathrm{e%
}^{\alpha\ln\xi}\rrvert \xi>0 \bigr) <\infty
\]
since the conditional distribution of $\ln\xi$ on $\xi>0$ belongs to the
class $\mathcal{S}(\alpha)$.

Our standing assumption is as follows:

%
\begin{assumption}
\label{Assumption A}Every convex combination of $F$ and $G$, namely $%
pF+(1-p)G$ for $0<p<1$, belongs to the class $\mathcal{R}_{-\alpha
}^{\ast}$.
\end{assumption}

Some interesting special cases of Assumption~\ref{Assumption A} include:
\begin{enumerate}[(b)]
\item[(a)] $F\in\mathcal{R}_{-\alpha}^{\ast}$ and $\overline{G}(x)=\mathrm{o}(
\overline{F}(x))$; or, symmetrically, $G\in\mathcal{R}_{-\alpha
}^{\ast}$
and $\overline{F}(x)=\mathrm{o}(\overline{G}(x))$.

\item[(b)] $F\in\mathcal{R}_{-\alpha}^{\ast}$, $G\in\mathcal
{R}_{-\alpha
}$, and $\overline{G}(x)=\mathrm{O}(\overline{F}(x))$; or, symmetrically,
$G\in
\mathcal{R}_{-\alpha}^{\ast}$, $F\in\mathcal{R}_{-\alpha}$, and $%
\overline{F}(x)=\mathrm{O}(\overline{G}(x))$.

\item[(c)] $F\in\mathcal{R}_{-\alpha}^{\ast}$, $G\in\mathcal
{R}_{-\alpha
}^{\ast}$, and the function $b(x)=\overline{F}(\mathrm
{e}^{x})/\overline{G}(%
\mathrm{e}^{x})$ is $O$-regularly varying (that is to say,
$b(xy)\asymp b(x)$
for every $y>0$).
\end{enumerate}
For (a) and (b), recall a fact that, if $V_{1}\in\mathcal{L}
(\alpha)$, $V_{2}\in\mathcal{L}(\alpha)$, and $\overline
{V_{1}}(x)\asymp
\overline{V_{2}}(x)$, then $V_{1}\in\mathcal{S}(\alpha)$ and
$V_{2}\in
\mathcal{S}(\alpha)$ are equivalent; see Theorem~2.1(a) of Kl\"{u}ppelberg
\cite{Klu88} and the sentences before it. This fact can be restated as that,
if $%
U_{1}\in\mathcal{R}_{-\alpha}$, $U_{2}\in\mathcal{R}_{-\alpha}$,
and $%
\overline{U_{1}}(x)\asymp\overline{U_{2}}(x)$, then $U_{1}\in
\mathcal
{R}%
_{-\alpha}^{\ast}$ and $U_{2}\in\mathcal{R}_{-\alpha}^{\ast}$ are
equivalent. By this fact, the verifications of (a) and (b) are
straightforward. For (c), by Theorem~2.0.8 of Bingham \textit{et al.} \cite{BinGolTeu87}, the
relation $b(xy)\asymp b(x)$ holds uniformly on every compact $y$-set of
$%
(0,\infty)$. Then the verification can be done by using Theorems 3.4 and
3.5 of Cline \cite{Cli87}.

\subsection{The main result}

In this subsection, we assume that $\{X,X_{1},X_{2},\ldots\}$ and $%
\{Y,Y_{1},Y_{2},\ldots\}$ are two independent sequences of i.i.d. random
variables with $X$ distributed as $F$ on $\mathbb{R}$ and $Y$ as $G$ on
$%
(0,\infty)$. Under Assumption~\ref{Assumption A}, by Lemma~\ref{C2} below
(with $n=2$), we have
\[
\Pr ( XY>x ) =\Pr ( X_{+}Y>x ) \sim\mathrm {E}Y^{\alpha}%
\overline{F}(x)+\mathrm{E}X_{+}^{\alpha}\overline{G}(x).
\]
Note that both $\mathrm{E}Y^{\alpha}$ and $\mathrm{E}X_{+}^{\alpha
}$ are
finite under Assumption~\ref{Assumption A}. The moments of $Y$ will appear
frequently in the paper, so we introduce a shorthand $\mu_{\alpha
}=\mathrm{%
E}Y^{\alpha}$ for $\alpha\geq0$ to help with the presentation. Starting
with this asymptotic formula and proceeding with induction, we shall
show in
our main result that the relations
%
%
\begin{equation}
\label{r1} \Pr ( M_{n}>x ) \sim A_{n}
\overline{F}(x)+B_{n}\overline{G}(x)
\end{equation}
and
%
%
\begin{equation}
\label{r2} \Pr ( S_{n}>x ) \sim A_{n}
\overline{F}(x)+C_{n}\overline{G}(x)
\end{equation}
hold for every $n\in\mathbb{N}$, where the coefficients $A_{n}$, $B_{n}$,
and $C_{n}$ are given by
\[
A_{n}=\sum_{i=1}^{n}
\mu_{\alpha}^{i},\qquad B_{n}=\sum
_{i=1}^{n}\mu _{\alpha}^{i-2}
\mathrm{E}M_{n-i+1}^{\alpha}, \qquad C_{n}=\sum
_{i=1}^{n}\mu _{\alpha}^{i-2}
\mathrm{E}S_{n-i+1,+}^{\alpha}.
\]

Furthermore, we shall seek to extend relations (\ref{r1}) and (\ref
{r2}) to $%
n=\infty$. For this purpose, it is natural to assume $\mu_{\alpha}<1$
(which excludes the case $\alpha=0$) to guarantee the finiteness of the
constants $A_{\infty}$, $B_{\infty}$, and $C_{\infty}$. Note in passing
that $\mu_{\alpha}<1$ implies $-\infty\leq\mathrm{E}\ln Y<0$,
which is
an aforementioned requirement for $S_{\infty}$ and $M_{\infty}$ to be a.s.
finite. Straightforwardly,
\[
A_{\infty}=\frac{\mu_{\alpha}}{1-\mu_{\alpha}}<\infty.
\]
It is easy to see that
%
%
\begin{equation}
\label{bound of moment} \mathrm{E}M_{\infty}^{\alpha}\leq\mathrm{E} \Biggl(
\sum_{i=1}^{\infty
}X_{i,+}\prod
_{j=1}^{i}Y_{j} \Biggr) ^{\alpha}<
\infty.
\end{equation}
Actually, when $0<\alpha\leq1$ we use the elementary inequality
$ (
\sum_{i=1}^{\infty}x_{i} ) ^{\alpha}\leq\sum_{i=1}^{\infty
}x_{i}^{\alpha}$ for any nonnegative sequence $\{x_{1},x_{2},\ldots\}$,
and when $\alpha>1$ we use Minkowski's inequality. In order for
$S_{\infty
} $ to be a.s. finite, we need another technical condition $\mathrm
{E}\ln
 ( X_{-}\vee1 ) <\infty$. The finiteness of $\mathrm
{E}S_{\infty
,+}^{\alpha}$ can be verified similarly to (\ref{bound of moment}).
Applying the dominated convergence theorem to the expressions for $B_{n}$
and $C_{n}$, we obtain
%
%
\begin{equation}
\label{BC} B_{\infty}=\frac{\mathrm{E}M_{\infty}^{\alpha}}{\mu_{\alpha
} (
1-\mu_{\alpha} ) }<\infty, \qquad C_{\infty}=
\frac{\mathrm{E}
S_{\infty,+}^{\alpha}}{\mu_{\alpha} ( 1-\mu_{\alpha} )
}%
<\infty.
\end{equation}

Now we are ready to state our main result, whose proof is postponed to
Sections~\ref{proof of (a)}--\ref{proof of (c)}.

%
\begin{theorem}
\label{No1}Let $\{X,X_{1},X_{2},\ldots\}$ and $\{Y,Y_{1},Y_{2},\ldots
\}$
be two independent sequences of i.i.d. random variables with $X$ distributed
as $F$ on $\mathbb{R}$ and $Y$ as $G$ on $(0,\infty)$. Under
Assumption~\ref{Assumption A}, we have the following:
\begin{enumerate}[(b)]
\item[(a)]Relations (\ref{r1}) and (\ref{r2}) hold for every $n\in\mathbb{N}$;

\item[(b)]If $\mu_{\alpha}<1$, then relation (\ref{r1}) holds for
$n=\infty$;

\item[(c)]If $\mu_{\alpha}<1$ and $\mathrm{E}\ln ( X_{-}\vee
1 )
<\infty$, then relation (\ref{r2}) holds for $n=\infty$.
\end{enumerate}
\end{theorem}

As we pointed out before, Theorem~\ref{No1} does not require a dominating
relationship between $\overline{F}$ and $\overline{G}$. Even in assertions
(b) and (c) where $\mu_{\alpha}<1$ is assumed, there is not
necessarily a
dominating relationship between $\overline{F}$ and $\overline{G}$, though
the conditions on $\overline{F}$ and $\overline{G}$ become not exactly
symmetric any more. Additionally, Theorems 5.2(3) and 6.1 of
Tang and Tsitsiashvili \cite{TanTsi03} are two special cases of
our Theorem~\ref{No1}(a)
with $%
\overline{G}(x)=\mathrm{o}(\overline{F}(x))$ and $\overline{F}(x)=\mathrm{o}(\overline
{G}(x))$, respectively.\vadjust{\goodbreak}

Since the famous work of Kesten \cite{Kes73}, the
tail probabilities of
$S_{\infty
}$ and $M_{\infty}$ have been extensively investigated, mainly in the
framework of random difference equations and most under so-called Cram\'
{e}r's condition that $\mu_{\alpha}=1$ holds for some $\alpha>0$.
Traditional random difference equations appearing in the literature are
often different from ones such as (\ref{M1}) and (\ref{T recursive
equation}) below associated to our model. Nevertheless, under our standard
assumptions on $\{X_{1},X_{2},\ldots\}$ and $\{Y_{1},Y_{2},\ldots\}$,
these subtle differences are not essential and the existing results can
easily be transformed to our framework. We omit such details here.
Corresponding to our model,  Kesten's work \cite{Kes73}
shows an asymptotic formula
of the form $Cx^{-\alpha}$ assuming, among others, that $Y$ fulfills
Cram\'{e}r's condition and $X$ fulfills a certain integrability condition
involving $Y$. Kesten's constant $C$, though positive, is generally unknown.
See Enriquez \textit{et al.} \cite{EnrSabZin09} for a probabilistic representation for this
constant. Goldie \cite{Gol91} studied the same
problem but in a broader scope and
he simplified Kesten's argument. Note that Cram\'{e}r's condition is
essentially used in these works. Among few works on this topic beyond
Cram\'{e}r's condition we mention Grey \cite{Gre94}
and Goldie and Gr\"{u}bel \cite{GolGru96}.
For the case where $F\in\mathcal{R}_{-\alpha}$ for some $\alpha>0$,
$\mu
_{\alpha+\varepsilon}<\infty$ for some $\varepsilon>0$, and $\mu
_{\alpha}<1$, indicating that the insurance risk dominates the financial
risk, Grey's work \cite{Gre94} shows a precise
asymptotic formula similar to ours.
Goldie and Gr\"{u}bel \cite{GolGru96} interpreted the study in terms of perpetuities
in insurance and finance and they derived some rough asymptotic formulas.
Corresponding to our model, their results show that $S_{\infty}$
exhibits a
light tail if $X$ is light tailed and $\Pr(Y\leq1)=1$, while
$S_{\infty}$
must exhibit a heavy tail once $\Pr(Y>1)>0$, regardless of the tail
behavior of $X$, all being consistent with the consensus on this topic that
risky investments are dangerous. We also refer the reader to
Hult and Samorodnitsky \cite{HulSam08}, Collamore
\cite{Col09}, Blanchet and Sigman \cite{BlaSig11}, and
Hitczenko and Weso\l owski \cite{HitWes11} for recent interesting developments
on the
topic.

In contrast to these existing results, we do not require Cram\'{e}r's
condition or a dominating relationship between $\overline{F}$ and
$\overline{%
G}$ in Theorem~\ref{No1}(b), (c). The coefficients $B_{\infty}$ and $%
C_{\infty}$ appearing in our formulas, though still generally unknown,
assume transparent structures as given in (\ref{BC}), which enable one to
easily conduct numerical estimates.

The condition $\mu_{\alpha}<1$ in Theorem~\ref{No1}(b), (c) is made mainly
to ensure the finiteness of $B_{\infty}$ and~$C_{\infty}$. However, it
excludes some apparently simpler cases such as $G\in\mathcal
{R}_{0}^{\ast}$
and classical random walks (corresponding to $\Pr(Y=1)=1$). The tail
behavior of the maximum of a random walk with negative drift, especially
with heavy-tailed increments, has been systematically investigated by many
people; see, for example, Feller \cite{Fel71},
Veraverbeke \cite{Ver77}, Korshunov \cite{Kor97},
Borovkov \cite{Bor03}, Denisov \textit{et al.} \cite{DenFosKor04}, and
Foss \textit{et al.} \cite{FosKorZac11}, among many
others. The study of random walks hints that the tail probabilities of $
S_{\infty}$ and $M_{\infty}$ behave essentially differently between the
cases $\mu_{\alpha}<1$ and $\mu_{\alpha}=1$. Actually, if $\mu
_{\alpha
}=1$, then all of $A_{n}$, $B_{n}$, and $C_{n}$ diverge to $\infty$ as
$%
n\rightarrow\infty$, and Theorem~\ref{No1} leads to
\[
\lim_{x\rightarrow\infty}\frac{\Pr ( S_{\infty}>x )
}{\overline{%
F}(x)+\overline{G}(x)}=\lim_{x\rightarrow\infty}
\frac{\Pr (
M_{\infty
}>x ) }{\overline{F}(x)+\overline{G}(x)}=\infty.
\]
This fails to give precise asymptotic formulas for $\Pr (
S_{\infty
}>x ) $ and $\Pr ( M_{\infty}>x ) $, though still consistent
with Kesten and Goldie's formula $Cx^{-\alpha}$ since $\overline{F}(x)+
\overline{G}(x)=\mathrm{o}(x^{-\alpha})$. For this case, intriguing questions
include how to capture the precise asymptotics other than Kesten and
Goldie's for $\Pr ( S_{\infty}>x ) $ and $\Pr (
M_{\infty
}>x ) $ and how to connect the asymptotics for $\Pr (
M_{n}>x ) $ and $\Pr ( S_{n}>x ) $ as $x\wedge
n\rightarrow
\infty$ to Kesten and Goldie's formula $Cx^{-\alpha}$. The approach
developed in the present paper seems not efficient to give a satisfactory
answer to either of these questions.

Admittedly, the standard complete independence assumptions on the two
sequences $\{X_{1},\allowbreak  X_{2},\ldots\}$ and $\{Y_{1},Y_{2},\ldots\}$, though
often appearing in the literature, are not of practical relevance. However,
Theorem~\ref{No1} offers new insights into the tail probabilities of the
sums in (\ref{S}) and their maxima in (\ref{M}), revealing the interplay
between the insurance and financial risks. Furthermore, extensions that
incorporate various dependence structures into the model are expected and
usually without much difficulty. We show in the next subsection a simple
example for such extensions.

\subsection{An extension}

As done by Chen \cite{Che11}, in this subsection
we assume that $%
\{(X,Y),(X_{1},Y_{1}),(X_{2},Y_{2}),\ldots\}$ is a sequence of i.i.d.
random pairs with $(X,Y)$ following a Farlie--Gumbel--Morgenstern (FGM)
distribution
%
%
\begin{equation}
\label{FGM} \pi(x,y)=F(x)G(y) \bigl( 1+\theta\overline{F}(x)\overline {G}(y)
\bigr) ,\qquad \theta\in{}[-1,1], x\in\mathbb{R}, y>0,
\end{equation}
where $F$ on $\mathbb{R}$ and $G$ on $(0,\infty)$ are two marginal
distributions. In view of the decomposition
%
%
\begin{equation}
\label{Pi} \pi=(1+\theta)FG-\theta F^{2}G-\theta
FG^{2}+\theta F^{2}G^{2},
\end{equation}
the FGM structure can easily be dissolved. Hereafter, for a random
variable $%
\xi$ and its i.i.d. copies $\xi_{1}$ and $\xi_{2}$, denote by
$\check
{\xi}
$ a random variable identically distributed as $\xi_{1}\vee\xi_{2}$ and
independent of all other sources of randomness. Under
Assumption~\ref{Assumption A}, by (\ref{Pi}) and Lemma~\ref{C2}
below, we can
conduct an
induction procedure to obtain
%
%
\begin{equation}
\label{r3} \Pr ( M_{n}>x ) \sim A_{n}^{\prime}
\overline {F}(x)+B_{n}^{\prime}%
\overline{G}(x)
\end{equation}
and
%
%
\begin{equation}
\label{r4} \Pr ( S_{n}>x ) \sim A_{n}^{\prime}
\overline {F}(x)+C_{n}^{\prime}%
\overline{G}(x)
\end{equation}
for every $n\in\mathbb{N}$, where
\begin{eqnarray*}
A_{n}^{\prime}& =& \bigl( ( 1-\theta ) \mu_{\alpha
}+\theta
\mathrm{E}\check{Y}^{\alpha} \bigr) \sum_{i=1}^{n}
\mu_{\alpha
}^{i-1},
\\
B_{n}^{\prime}& =&\sum_{i=1}^{n}
\mu_{\alpha}^{i-1} \bigl( ( 1-\theta ) \mathrm{E} (
M_{n-i}+X_{n-i+1} ) _{+}^{\alpha
}+\theta
\mathrm{E} ( M_{n-i}+\check{X}_{n-i+1} ) _{+}^{\alpha
}
\bigr),
\\
C_{n}^{\prime}& =&\sum_{i=1}^{n}
\mu_{\alpha}^{i-1} \bigl( ( 1-\theta ) \mathrm{E} (
S_{n-i}+X_{n-i+1} ) _{+}^{\alpha
}+\theta
\mathrm{E} ( S_{n-i}+\check{X}_{n-i+1} ) _{+}^{\alpha
}
\bigr).
\end{eqnarray*}
Additionally, under the conditions of Theorem~\ref{No1}(b), (c), letting $
n\rightarrow\infty$ leads to
\begin{eqnarray*}
A_{\infty}^{\prime}& =&\frac{1}{1-\mu_{\alpha}} \bigl( ( 1-\theta )
\mu_{\alpha}+\theta\mathrm{E}\check{Y}^{\alpha} \bigr) ,
\\
B_{\infty}^{\prime}& =&\frac{1}{1-\mu_{\alpha}} \bigl( ( 1-\theta )
\mathrm{E} ( M_{\infty}+X ) _{+}^{\alpha}+\theta
\mathrm{E%
} ( M_{\infty}+\check{X} ) _{+}^{\alpha}
\bigr),
\\
C_{\infty}^{\prime}& =&\frac{1}{1-\mu_{\alpha}} \bigl( ( 1-\theta )
\mathrm{E} ( S_{\infty}+X ) _{+}^{\alpha}+\theta
\mathrm{E%
} ( S_{\infty}+\check{X} ) _{+}^{\alpha}
\bigr),
\end{eqnarray*}
where $X$ and $\check{X}$ are independent of $M_{\infty}$ and
$S_{\infty}$. It is easy to verify the finiteness of $B_{\infty}^{\prime}$ and $%
C_{\infty}^{\prime}$.

We summarize the analysis above into the following corollary and will
show a
sketch of its proof in Section~\ref{proof of cor}.

%
\begin{corollary}
\label{No2}Let $\{(X,Y),(X_{1},Y_{1}),(X_{2},Y_{2}),\ldots\}$ be a sequence
of i.i.d. random pairs with common FGM distribution (\ref{FGM}). Under
Assumption~\ref{Assumption A}, we have the following:
\begin{enumerate}[(b)]
\item[(a)]Relations (\ref{r3}) and (\ref{r4}) hold for every $n\in\mathbb{N}$;

\item[(b)]If $\mu_{\alpha}<1$, then relation (\ref{r3}) holds for
$n=\infty$;

\item[(c)]If $\mu_{\alpha}<1$ and $\mathrm{E}\ln ( X_{-}\vee
1 )
<\infty$, then relation (\ref{r4}) holds for $n=\infty$.
\end{enumerate}
\end{corollary}

As a sanity check, letting $\theta=0$, the results in Corollary~\ref{No2}
coincide with those in Theorem~\ref{No1}.

\section{Lemmas}\label{sec3}

In this section, we prepare a series of lemmas, some of which are
interesting in their own right. We first recall some well-known properties
of distributions of regular variation and convolution equivalence. If
$U\in
\mathcal{R}_{-\alpha}$ for some $0\leq\alpha<\infty$, then for
every $%
\varepsilon>0$ and every $b>1$ there is some constant $x_{0}>0$ such that
Potter's bounds
%
%
\begin{equation}
\label{Potter's} \frac{1}{b} \bigl( y^{-\alpha-\varepsilon}\wedge y^{-\alpha
+\varepsilon
}
\bigr) \leq\frac{\overline{U}(xy)}{\overline{U}(x)}\leq b \bigl( y^{-\alpha-\varepsilon}\vee y^{-\alpha+\varepsilon}
\bigr)
\end{equation}
hold whenever $x\geq x_{0}$ and $xy\geq x_{0}$; see Theorem~1.5.6(iii) of
Bingham \textit{et al.} \cite{BinGolTeu87}. Since $U\in\mathcal{R}_{-\alpha}$ if and only
if $V$
defined by (\ref{U--V}) belongs to $\mathcal{L}(\alpha)$, Potter's bounds
above can easily be restated in terms of a distribution $V\in\mathcal
{L}%
(\alpha)$ as that, for every $\varepsilon>0$ and every $b>1$ there is some
constant $x_{0}>0$ such that the inequalities
%
%
\begin{equation}
\label{Potter's for L} \frac{1}{b} \bigl( \mathrm{e}^{- ( \alpha+\varepsilon )
y}\wedge
\mathrm{e}^{- ( \alpha-\varepsilon ) y} \bigr) \leq \frac{%
\overline{V}(x+y)}{\overline{V}(x)}\leq b \bigl(
\mathrm{e}^{- (
\alpha
+\varepsilon ) y}\vee\mathrm{e}^{- ( \alpha-\varepsilon
 )
y} \bigr)
\end{equation}
hold whenever $x\geq x_{0}$ and $x+y\geq x_{0}$. By Lemma~5.2 of
Pakes \cite{Pak04}, if $V\in\mathcal{S}(\alpha
)$ then it holds for every $n\in
\mathbb{%
N}$ that
%
%
\begin{equation}
\label{Pakes} \lim_{x\rightarrow\infty}\frac{\overline{V^{n\ast}}(x)}{\overline
{V}(x)}%
=n \bigl(
\hat{V}(\alpha) \bigr) ^{n-1}.
\end{equation}

The first lemma below describes an elementary property of convolution
equivalence.

%
\begin{lemma}
\label{L1}Let $\eta_{1}, \ldots, \eta_{n}$ be $n\geq2$ i.i.d.
real-valued random variables with common distribution $V\in\mathcal{S}
(\alpha)$ for some $\alpha\geq0$. Then
\[
\lim_{c\rightarrow\infty}\lim_{x\rightarrow\infty}\frac{\Pr
 (
\sum_{i=1}^{n}\eta_{i}>x,\eta_{1}>c,\eta_{2}>c ) }{\overline
{V}(x)}%
=0.
\]
\end{lemma}

\begin{pf}
For every $x\geq0$ and $c\geq0$, write
%
%
\begin{eqnarray}
\label{123} && \Pr \Biggl( \sum_{i=1}^{n}
\eta_{i}>x,\eta_{1}>c,\eta_{2}>c \Biggr)
\nonumber
\\
&&\quad =\Pr \Biggl( \sum_{i=1}^{n}
\eta_{i}>x \Biggr) -2\Pr \Biggl( \sum_{i=1}^{n}
\eta_{i}>x,\eta_{1}\leq c \Biggr) +\Pr \Biggl( \sum
_{i=1}^{n}\eta_{i}>x,\eta_{1}\leq
c,\eta_{2}\leq c \Biggr)\quad
\\
&&\quad =I_{1}(x)-2I_{2}(x,c)+I_{3}(x,c).
\nonumber
\end{eqnarray}
By relation (\ref{Pakes}), we have
\[
\lim_{x\rightarrow\infty}\frac{I_{1}(x)}{\overline{V}(x)}=n \bigl( \hat {V}%
(
\alpha) \bigr) ^{n-1}
\]
and
\begin{eqnarray*}
\lim_{x\rightarrow\infty}\frac{I_{2}(x,c)}{\overline{V}(x)}& =&\lim_{x\rightarrow\infty}
\int_{-\infty}^{c}\frac{\Pr (
\sum_{i=1}^{n-1}\eta_{i}>x-y ) }{\overline{V}(x-y)}\frac
{\overline
{V}%
(x-y)}{\overline{V}(x)}V(
\mathrm{d}y)
\\
& =&(n-1) \bigl( \hat{V}(\alpha) \bigr) ^{n-2}\int_{-\infty
}^{c}
\mathrm {e}%
^{\alpha y}V(\mathrm{d}y),
\end{eqnarray*}
where in the last step we used $V\in\mathcal{L}(\alpha)$ and the dominated
convergence theorem. Similarly,
\begin{eqnarray*}
\lim_{x\rightarrow\infty}\frac{I_{3}(x,c)}{\overline{V}(x)}& =&\lim_{x\rightarrow\infty}
\int_{-\infty}^{c}\int_{-\infty
}^{c}
\frac
{\Pr
 ( \sum_{i=1}^{n-2}\eta_{i}>x-y_{1}-y_{2} ) }{\overline
{V}(x)}V(%
\mathrm{d}y_{1})V(\mathrm{d}y_{2})
\\
& =&(n-2) \bigl( \hat{V}(\alpha) \bigr) ^{n-3} \biggl( \int
_{-\infty
}^{c}%
\mathrm{e}^{\alpha y}V(
\mathrm{d}y) \biggr) ^{2}.
\end{eqnarray*}
Plugging these limits into (\ref{123}) yields the desired result.
\end{pf}

Hereafter, for $n\geq2$ distributions $V_{1}, \ldots, V_{n}$,
denote by $%
V_{\mathbf{p}}=\sum_{i=1}^{n}p_{i}V_{i}$ a convex combination of $V_{1},
\ldots, V_{n}$, where $\mathbf{p}\in\Delta=\{(p_{1},\ldots
,p_{n})\in
(0,1)^{n}\dvtx \sum_{i=1}^{n}p_{i}=1\}$.

%
\begin{lemma}
\label{new1}Let $V_{1}, \ldots, V_{n}$ be $n\geq2$ distributions
and let
$\alpha\geq0$. The following are equivalent:
\begin{enumerate}[(b)]
\item[(a)] $V_{\mathbf{p}}\in\mathcal{S}(\alpha)$ for every
$\mathbf
{p}\in
\Delta$;

\item[(b)] $V_{\mathbf{p}}\in\mathcal{S}(\alpha)$ for some
$\mathbf
{p}\in
\Delta$ and the relation
%
%
\begin{equation}
\label{assumption} \overline{V_{i}}(x-y)-\mathrm{e}^{\alpha y}\overline
{V_{i}}(x)=\mathrm{o} \Biggl( \sum_{j=1}^{n}
\overline{V_{j}}(x) \Biggr)
\end{equation}
holds for every $y\in\mathbb{R}$ and every $i=1,\ldots,n$.
\end{enumerate}
\end{lemma}

\begin{pf}
First prove that (b) implies (a). Denote by $\mathbf{p}^{\ast}$ this
specific element in $\Delta$ such that $V_{\mathbf{p}^{\ast}}\in
\mathcal{S%
}(\alpha)$. For every $\mathbf{p}\in\Delta$, it is easy to see that $
\overline{V_{\mathbf{p}}}(x)\asymp\sum_{j=1}^{n}\overline
{V_{j}}(x)\asymp
\overline{V_{\mathbf{p}^{\ast}}}(x)$ and that $V_{\mathbf{p}}\in
\mathcal{L}%
(\alpha)$ by (\ref{assumption}). Thus, $V_{\mathbf{p}}\in\mathcal
{S}%
(\alpha)$ follows from the closure of the class $\mathcal{S}(\alpha)$
under weak equivalence as mentioned in the last paragraph of Section~\ref{A}.

For the other implication, we only need to use (a) to verify
(\ref{assumption}). For arbitrarily fixed $0<\varepsilon<1$ and every $%
i=1,\ldots,n$, each of the sums $\overline{V_{i}}(x)+\varepsilon
\sum_{j=1,j\neq i}^{n}\overline{V_{j}}(x)$ and $\sum_{j=1}^{n}\overline
{V_{j}%
}(x)$ is proportional to a convolution-equivalent tail. Thus,
\begin{eqnarray*}
&&\bigl\llvert \overline{V_{i}}(x-y)-\mathrm{e}^{\alpha y}\overline
{V_{i}}%
(x)\bigr\rrvert
\\
&&\quad  \leq \Biggl\llvert \bigl( \overline
{V_{i}}(x-y)-\mathrm{e} ^{\alpha y}\overline{V_{i}}(x)
\bigr) +\varepsilon\sum_{j=1,j\neq
i}^{n} \bigl(
\overline{V_{j}}(x-y)-\mathrm{e}^{\alpha y}\overline{V_{j}}
(x) \bigr) \Biggr\rrvert
\\
&& \qquad {} +\varepsilon\sum_{j=1,j\neq i}^{n}\bigl\llvert
\overline {V_{j}}(x-y)-%
\mathrm{e}^{\alpha y}
\overline{V_{j}}(x)\bigr\rrvert
\\
&&\quad  \leq \Biggl\llvert \Biggl( \overline{V_{i}}(x-y)+\varepsilon\sum
_{j=1,j\neq
i}^{n}\overline{V_{j}}(x-y)
\Biggr) -\mathrm{e}^{\alpha y} \Biggl( \overline{%
V_{i}}(x)+
\varepsilon\sum_{j=1,j\neq i}^{n}
\overline{V_{j}}(x) \Biggr) \Biggr\rrvert
\\
&& \qquad {} +\varepsilon\sum_{j=1}^{n}
\overline{V_{j}}(x-y)+\varepsilon \mathrm{%
e}^{\alpha y}
\sum_{j=1}^{n}\overline{V_{j}}(x)
\\
&&\quad  = \mathrm{o} ( 1 ) \Biggl( \overline{V_{i}}(x)+\varepsilon\sum
_{j=1,j\neq
i}^{n}\overline{V_{j}}(x) \Biggr) +2
\varepsilon \bigl( \mathrm {e}^{\alpha
y}+\mathrm{o}(1) \bigr) \sum
_{j=1}^{n}\overline{V_{j}}(x).
\end{eqnarray*}
By the arbitrariness of $\varepsilon$, relation (\ref{assumption})
follows.
\end{pf}

The following lemma shows the usefulness of convolution equivalence in
dealing with the tail probability of the sum of independent random
variables. Note\vadjust{\goodbreak} that the lemma does not require any dominating relationship
among the individual tails. Additionally, in view of Lemma~\ref{new1},
letting $\alpha=0$ in Lemma~\ref{L2} retrieves Theorem~1 of
Li and Tang \cite{LiTan10}.

%
\begin{lemma}
\label{L2}Let $V_{1}, \ldots, V_{n}$ be $n\geq2$ distributions on $
\mathbb{R}$ and let $\alpha\geq0$. If $V_{\mathbf{p}}\in\mathcal
{S}%
(\alpha)$ for every $\mathbf{p}\in\Delta$, then $V_{1}\ast\cdots
\ast
V_{n}\in\mathcal{S}(\alpha)$ and
%
%
\begin{equation}
\label{for n} \overline{V_{1}\ast\cdots\ast V_{n}}(x)\sim
\sum_{i=1}^{n} \Biggl( \prod
_{j=1,j\neq i}^{n}\hat{V}_{j}(\alpha) \Biggr)
\overline{V_{i}}(x).
\end{equation}
\end{lemma}

\begin{pf}
Clearly, we only need to prove relation (\ref{for n}). Introduce $n$
independent random variables $\eta_{1}, \ldots, \eta_{n}$ with
distributions $V_{1}, \ldots, V_{n}$, respectively. For every
$x\geq0$
and $0\leq c\leq x/n$,
\[
\overline{V_{1}\ast\cdots\ast V_{n}}(x)=\Pr \Biggl( \sum
_{i=1}^{n}\eta _{i}>x,\bigcup
_{j=1}^{n} ( \eta_{j}>c ) \Biggr).
\]
According to whether or not there is exactly only one $(\eta_{j}>c)$
occurring in the union, we split the probability on the right-hand side into
two parts as
%
%
\begin{equation}
\label{c} \overline{V_{1}\ast\cdots\ast V_{n}}(x)=I_{1}(x,c)+I_{2}(x,c).
\end{equation}
First we deal with $I_{1}(x,c)$. For a real vector $\mathbf
{y}=(y_{1},\ldots
,y_{n-1})^{\prime}$, write $\Sigma=\sum_{k=1}^{n-1}y_{k}$, and for
each $%
j=1,\ldots,n$, write
\[
\Biggl( \prod_{k=1,k\neq j}^{n}
\mathrm{d}V_{k} \Biggr) (\mathbf{y})=V_{1}(
\mathrm{d}y_{1})\cdots V_{j-1}(\mathrm{d}y_{j-1})V_{j+1}(
\mathrm{d}%
y_{j})\cdots V_{n}(\mathrm{d}y_{n-1}).
\]
We have
\begin{eqnarray*}
I_{1}(x,c)& =&\sum_{j=1}^{n}\Pr
\Biggl( \sum_{i=1}^{n}\eta_{i}>x,
\eta _{j}>c,\bigcap_{k=1,k\neq j}^{n} (
\eta_{k}\leq c ) \Biggr)
\\
& =&\sum_{j=1}^{n}\int_{-\infty}^{c}
\cdots\int_{-\infty
}^{c}\overline {V_{j}%
}
( x-\Sigma ) \Biggl( \prod_{k=1,k\neq j}^{n}
\mathrm{d}%
V_{k} \Biggr) (\mathbf{y})
\\
& =&\int_{-\infty}^{c}\cdots\int_{-\infty}^{c}
\Biggl( \sum_{j=1}^{n}%
\overline{V_{j}} ( x-\Sigma ) \Biggr) \Biggl( \sum
_{h=1}^{n} \Biggl( \prod_{k=1,k\neq h}^{n}
\mathrm{d}V_{k} \Biggr) (\mathbf{y}) \Biggr)
\\
&&{} -\sum_{j=1}^{n}\sum
_{h=1,h\neq j}^{n}\int_{-\infty
}^{c}
\cdots \int_{-\infty}^{c}\overline{V_{j}} (
x-\Sigma ) \Biggl( \prod_{k=1,k\neq h}^{n}
\mathrm{d}V_{k} \Biggr) (\mathbf{y}).
\end{eqnarray*}
Since $\sum_{j=1}^{n}\overline{V_{j}} ( x ) $ is proportional
to a
convolution-equivalent tail, by the dominated convergence theorem,
\begin{eqnarray*}
I_{1}(x,c)& \sim& \Biggl( \sum_{j=1}^{n}
\overline{V_{j}} ( x ) \Biggr) \int_{-\infty}^{c}
\cdots\int_{-\infty}^{c}\mathrm {e}^{\alpha
\Sigma}
\Biggl( \sum_{h=1}^{n} \Biggl( \prod
_{k=1,k\neq h}^{n}\mathrm {d}%
V_{k}
\Biggr) (\mathbf{y}) \Biggr)
\\
&&{} -\sum_{j=1}^{n}\sum
_{h=1,h\neq j}^{n}\int_{-\infty
}^{c}
\cdots \int_{-\infty}^{c}\overline{V_{j}} (
x-\Sigma ) \Biggl( \prod_{k=1,k\neq h}^{n}
\mathrm{d}V_{k} \Biggr) (\mathbf{y})
\\
& =&\sum_{j=1}^{n}\overline{V_{j}}
( x ) \int_{-\infty
}^{c}\cdots \int
_{-\infty}^{c}\mathrm{e}^{\alpha\Sigma} \Biggl( \prod
_{k=1,k\neq
j}^{n}%
\mathrm{d}V_{k}
\Biggr) (\mathbf{y})
\\
&&{} -\sum_{j=1}^{n}\sum
_{h=1,h\neq j}^{n}\int_{-\infty
}^{c}
\cdots \int_{-\infty}^{c} \bigl( \overline{V_{j}}
( x-\Sigma ) -\mathrm{e}%
^{\alpha\Sigma}\overline{V_{j}} ( x
) \bigr) \Biggl( \prod_{k=1,k\neq h}^{n}
\mathrm{d}V_{k} \Biggr) (\mathbf{y}).
\end{eqnarray*}
Hence, it follows from (\ref{assumption}) and the dominated convergence
theorem that
%
%
\begin{equation}
\label{a} \lim_{c\rightarrow\infty}\lim_{x\rightarrow\infty}
\frac{I_{1}(x,c)}{
\sum_{i=1}^{n} ( \prod_{j=1,j\neq i}^{n}\hat{V}_{j}(\alpha
) )
\overline{V_{i}}(x)}=1.
\end{equation}
Next we turn to $I_{2}(x,c)$. Write $\tilde{\eta}=\max\{\eta
_{1},\ldots
,\eta_{n}\}$, which has a convolution-equivalent tail proportional to $
\sum_{j=1}^{n}\overline{V_{j}} ( x ) $, and let $\tilde
{\eta}_{1},
\ldots, \tilde{\eta}_{n}$ be i.i.d. copies of $\tilde{\eta}$. Clearly,
\begin{eqnarray*}
I_{2}(x,c)& =&\Pr \Biggl( \sum_{i=1}^{n}
\eta_{i}>x,\bigcup_{1\leq
j<k\leq
n} (
\eta_{j}>c,\eta_{k}>c ) \Biggr)
\\
& \leq&\sum_{1\leq j<k\leq n}\Pr \Biggl( \sum
_{i=1}^{n}\tilde{\eta }_{i}>x,%
\tilde{\eta}_{j}>c,\tilde{\eta}_{k}>c \Biggr).
\end{eqnarray*}
Thus, by Lemma~\ref{L1},
%
%
\begin{eqnarray}
\label{b} && \lim_{c\rightarrow\infty}\limsup_{x\rightarrow\infty}
\frac
{I_{2}(x,c)}{%
\sum_{i=1}^{n} ( \prod_{j=1,j\neq i}^{n}\hat{V}_{j}(\alpha
) )
\overline{V_{i}}(x)}
\nonumber
\\
&&\quad \leq\lim_{c\rightarrow\infty}\lim_{x\rightarrow\infty}
\frac
{I_{2}(x,c)%
}{\sum_{j=1}^{n}\overline{V_{j}} ( x ) }\limsup_{x\rightarrow
\infty}\frac{\sum_{j=1}^{n}\overline{V_{j}} ( x ) }{%
\sum_{i=1}^{n} ( \prod_{j=1,j\neq i}^{n}\hat{V}_{j}(\alpha
) )
\overline{V_{i}}(x)}
\\
&&\quad =0.
\nonumber
\end{eqnarray}
Plugging (\ref{a}) and (\ref{b}) into (\ref{c}) yields the desired
result.
\end{pf}

Due to the connection between convolution equivalence and strongly regular
variation, we can restate Lemmas~\ref{new1} and~\ref{L2} in terms of
strongly regular variation. Actually, the next lemma shows an equivalent
condition for Assumption~\ref{Assumption A}.

%
\begin{lemma}
\label{new2}Let $U_{1}, \ldots, U_{n}$ be $n\geq2$ distributions
and let
$\alpha\geq0$. The following are equivalent:
\begin{enumerate}[(b)]
\item[(a)] $U_{\mathbf{p}}\in\mathcal{R}_{-\alpha}^{\ast}$ for every
$%
\mathbf{p}\in\Delta$;

\item[(b)] $U_{\mathbf{p}}\in\mathcal{R}_{-\alpha}^{\ast}$ for
some $
\mathbf{p}\in\Delta$ and the relation
\[
\overline{U_{i}}(x/y)-y^{\alpha}\overline{U_{i}}(x)=\mathrm{o}
\Biggl( \sum_{j=1}^{n}%
\overline{U_{j}}(x) \Biggr)
\]
holds for every $y>0$ and every $i=1,\ldots,n$.
\end{enumerate}
\end{lemma}

The lemma below expands the tail probability of the product of independent,
nonnegative, and strongly regular random variables, forming an analogue of
the well-known Breiman's theorem in a different situation. For Breiman's
theorem, see Breiman \cite{Bre65} and Cline and Samorodnitsky \cite{CliSam94}.

%
\begin{lemma}
\label{C2}Let $\xi_{1}, \ldots, \xi_{n}$ be $n\geq2$ independent
nonnegative random variables with distributions $U_{1}, \ldots, U_{n}$,
respectively, and let $\alpha\geq0$. If $U_{\mathbf{p}}\in\mathcal{R}
_{-\alpha}^{\ast}$ for every $\mathbf{p}\in\Delta$, then the
distribution of $\prod_{i=1}^{n}\xi_{i}$ belongs to the class
$\mathcal
{R}%
_{-\alpha}^{\ast}$ and
\[
\Pr \Biggl( \prod_{i=1}^{n}
\xi_{i}>x \Biggr) \sim\sum_{i=1}^{n}
\Biggl( \prod_{j=1,j\neq i}^{n}\mathrm{E}
\xi_{j}^{\alpha} \Biggr) \overline {U_{i}}%
(x).
\]
\end{lemma}

The next lemma shows Kesten's bound for convolution tails without the usual
requirement $\hat{V}(\alpha)\geq1$. It improves Lemma~5.3 of
Pakes \cite{Pak04}
for the case $0<\hat{V}(\alpha)<1$.

%
\begin{lemma}
\label{Kesten}Let $V$ be a distribution on $\mathbb{R}$. If $V\in
\mathcal{S}%
(\alpha)$ for some $\alpha\geq0$, then for every $\varepsilon>0$ there
is some constant $K>0$ such that the relation
\[
\overline{V^{n\ast}}(x)\leq K \bigl( \hat{V}(\alpha)+\varepsilon \bigr)
^{n}\overline{V}(x)
\]
holds for all $n\in\mathbb{N}$ and all $x\geq0$.
\end{lemma}

\begin{pf}
When $\hat{V}(\alpha)\geq1$, the assertion has been given in Lemma~5.3 of
Pakes \cite{Pak04}. Hence, we only need to
consider $\hat{V}(\alpha)<1$ (for
which $\alpha>0$ must hold). Let $\{\eta,\eta_{1},\eta_{2},\ldots\}
$ be
a sequence of i.i.d. random variables with common distribution $V$, and
set $%
c=-\alpha^{-1}\ln\hat{V}(\alpha)>0$. Clearly,
\[
\overline{V^{n\ast}}(x)=\Pr \Biggl( \sum_{i=1}^{n}
( \eta _{i}+c ) >x+nc \Biggr).
\]
Note that the distribution of $\eta+c$ still belongs to the class
$\mathcal{%
S}(\alpha)$ and $\mathrm{Ee}^{\alpha(\eta+c)}=1$. Hence, for every $
\delta>0$, by Lemma~5.3 of Pakes \cite{Pak04},
there is some constant $K_{1}>0$
such that, for all $n\in\mathbb{N}$ and all $x\geq0$,
%
%
\begin{equation}
\label{k1} \overline{V^{n\ast}}(x)\leq K_{1}(1+
\delta)^{n}\Pr ( \eta +c>x+nc ) =K_{1}(1+
\delta)^{n}\overline{V} \bigl( x+(n-1)c \bigr).
\end{equation}
By (\ref{Potter's for L}), there are some constants $K_{2}>0$ and $x_{0}>0$
such that, for all $n\in\mathbb{N}$ and all $x\geq x_{0}$,
%
%
\begin{equation}
\label{k2} \overline{V} \bigl( x+(n-1)c \bigr) \leq K_{2}
\mathrm{e}^{-(\alpha
-\delta
)(n-1)c}\overline{V} ( x ).
\end{equation}
Plugging (\ref{k2}) into (\ref{k1}) and noticing that $\mathrm
{e}^{-\alpha
c}=\hat{V}(\alpha)$, we have, for all $n\in\mathbb{N}$ and all
$x\geq
x_{0} $,
%
%
\begin{equation}
\label{k3} \overline{V^{n\ast}}(x)\leq K_{1}K_{2}
\mathrm{e}^{(\alpha-\delta
)c} \bigl( (1+\delta)\mathrm{e}^{c\delta}\hat{V}(
\alpha) \bigr) ^{n}\overline{V}(x).
\end{equation}
For $0\leq x<x_{0}$, we choose an integer $n_{0}\geq x_{0}/c$. Then,
for $%
0\leq x<x_{0}$ and $n>n_{0}$, using the same derivations as in (\ref
{k1})--(\ref{k3}), we obtain
%
%
\begin{eqnarray}
\label{k4} \overline{V^{n\ast}}(x)& \leq& K_{1}(1+
\delta)^{n}\overline{V} \bigl( x+n_{0}c+(n-n_{0}-1)c
\bigr)
\nonumber
\\
& \leq& K_{1}K_{2}\mathrm{e}^{(\alpha-\delta)(n_{0}+1)c} \bigl( (1+
\delta)%
\mathrm{e}^{c\delta}\hat{V}(\alpha) \bigr) ^{n}
\overline{V} ( x+n_{0}c )
\\
& \leq& K_{1}K_{2}\mathrm{e}^{(\alpha-\delta)(n_{0}+1)c} \bigl( (1+
\delta)%
\mathrm{e}^{c\delta}\hat{V}(\alpha) \bigr) ^{n}
\overline{V} ( x ) .
\nonumber
\end{eqnarray}
At last, for $0\leq x<x_{0}$ and $1\leq n\leq n_{0}$, it is obvious that
%
%
\begin{equation}
\label{k5} \overline{V^{n\ast}}(x)\leq1\leq\frac{ ( (1+\delta)\mathrm
{e}%
^{c\delta}\hat{V}(\alpha) ) ^{n}}{ ( (1+\delta)\mathrm
{e}%
^{c\delta}\hat{V}(\alpha) ) ^{n_{0}}\wedge1}
\frac{\overline
{V} (
x ) }{\overline{V} ( x_{0} ) }.
\end{equation}
A combination of (\ref{k3})--(\ref{k5}) gives that, for some constant $K>0$
and for all $n\in\mathbb{N}$ and all $x\geq0$,
\[
\overline{V^{n\ast}}(x)\leq K \bigl( (1+\delta)\mathrm{e}^{c\delta
}
\hat{V}%
(\alpha) \bigr) ^{n}\overline{V} ( x ).
\]
By setting $\delta$ to be small enough such that $(1+\delta)\mathrm{e}
^{c\delta}\hat{V}(\alpha)\leq\hat{V}(\alpha)+\varepsilon$, we complete
the proof.
\end{pf}

The following lemma will be crucial in proving Theorem~\ref{No1}(b), (c).

%
\begin{lemma}
\label{infinite}Let $\{X,X_{1},X_{2},\ldots\}$ be a sequence of
(arbitrarily dependent) random variables with common distribution $F$
on $%
\mathbb{R}$, let $\{Y,Y_{1},Y_{2},\ldots\}$ be another sequence of i.i.d.
random variables with common distribution $G$ on $[0,\infty)$, and let the
two sequences be mutually independent. Assume that there is some
distribution $U\in\mathcal{R}_{-\alpha}^{\ast}$ for $\alpha>0$
such that
\[
\overline{F}(x)+\overline{G}(x)=\mathrm{O}\bigl(\overline{U}(x)\bigr).
\]
Assume also that $\mu_{\alpha}<1$. Then
%
%
\begin{equation}
\label{bound} \lim_{n\rightarrow\infty}\limsup_{x\rightarrow\infty}
\frac
{1}{\overline{U%
}(x)}\Pr \Biggl( \sum_{i=n+1}^{\infty}X_{i}
\prod_{j=1}^{i}Y_{j}>x \Biggr) =0.
\end{equation}
\end{lemma}

\begin{pf}
Choose some large constant $K_{1}>0$ such that the inequality
$\overline
{F}%
(x)\vee\overline{G}(x)\leq K_{1}\overline{U}(x)$ holds for all $x\in
\mathbb{R}$, and then introduce a nonnegative random variable $X^{\ast}$
with a distribution
\[
F^{\ast}(x)= \bigl( 1-K_{1}\overline{U}(x) \bigr)
_{+},\qquad x\geq0.
\]
Clearly, $\overline{F}(x)\leq\overline{F^{\ast}}(x)\leq
K_{1}\overline
{U}%
(x)$ for all $x\geq0$ and $\overline{F^{\ast}}(x)=K_{1}\overline{U}(x)$
for all large $x$. The inequality $\overline{F}(x)\leq\overline
{F^{\ast}}%
(x)$ for all $x\geq0$ means that $X$ is stochastically not greater
than $%
X^{\ast}$, denoted by $X\leq_{\mathrm{st}}X^{\ast}$. Moreover, since
$%
U\in\mathcal{R}_{-\alpha}^{\ast}$, there is some large but fixed constant
$t>0$ such that $K_{1}\int_{t}^{\infty}z^{\alpha}U(\mathrm
{d}z)<1-\mu
_{\alpha}$. For this fixed $t$, define
\[
t_{0}=\inf\bigl\{s\geq t\dvtx K_{1}\overline{U}(s)\leq
\overline{G}(t)\bigr\},
\]
and then introduce another nonnegative random variable $Y^{\ast}$ with a
distribution
\[
G^{\ast}(x)=G(x)\mathbf{1}_{(0\leq x<t)}+G(t)\mathbf{1}_{(t\leq
x<t_{0})}+
\bigl( 1-K_{1}\overline{U}(x) \bigr) \mathbf{1}_{(x\geq t_{0})}.
\]
Clearly, $\mathrm{E} ( Y^{\ast} ) ^{\alpha}<1$,
$\overline{G}
(x)\leq\overline{G^{\ast}}(x)\leq K_{1}\overline{U}(x)$ for all
$x>0$, and
$\overline{G^{\ast}}(x)=K_{1}\overline{U}(x)$ for all \mbox{$x\geq t_{0}$}.
Thus, $%
Y\leq_{\mathrm{st}}Y^{\ast}$. Let $Y_{1}^{\ast}$, $Y_{2}^{\ast}$,
\ldots\ be i.i.d. copies of $Y^{\ast}$ independent of $X^{\ast}$.

Choose some $0<\varepsilon<\alpha\wedge ( 1-\mathrm{E} (
Y^{\ast
} ) ^{\alpha} ) $ such that $\mathrm{E} ( Y^{\ast
} )
^{\alpha-\varepsilon}<1$. By Lemma~\ref{Kesten}, there is some
constant $%
K_{2}>0$ such that, for all $i\in\mathbb{N}$ and all $x\geq1$,
%
%
\begin{equation}
\label{inf1} \Pr \Biggl( \prod_{j=1}^{i}Y_{j}^{\ast}>x
\Biggr) =\Pr \Biggl( \sum_{j=1}^{i}\ln
Y_{j}^{\ast}>\ln x \Biggr) \leq K_{2} \bigl(
\mathrm{E} \bigl( Y^{\ast} \bigr) ^{\alpha}+\varepsilon \bigr)
^{i}\overline {G^{\ast
}}(x).
\end{equation}
Noticeably, the derivation in (\ref{inf1}) tacitly requires that $%
Y_{1}^{\ast}, \ldots, Y_{j}^{\ast}$ are positive. Nevertheless, in case
$G^{\ast}$ assigns a mass at $0$, the upper bound in (\ref{inf1}) is still
correct and can easily be verified by conditioning on $\bigcap_{j=1}^{i}
 ( Y_{j}^{\ast}>0 ) $. By Lemma~\ref{C2},
%
%
\begin{equation}
\label{inf21} \Pr \bigl( X^{\ast}Y^{\ast}>x \bigr) \sim
K_{1} \bigl( \mathrm {E} \bigl( X^{\ast} \bigr)
^{\alpha}+\mathrm{E} \bigl( Y^{\ast} \bigr) ^{\alpha
} \bigr)
\overline{U}(x).
\end{equation}
Moreover, by (\ref{Potter's}), there is some constant $x_{0}>0$ such that,
for all $x>x_{0}$ and $xy>x_{0}$,
%
%
\begin{equation}
\label{inf22} \overline{U}(xy)\leq(1+\varepsilon) \bigl( y^{-\alpha-\varepsilon
}\vee
y^{-\alpha+\varepsilon} \bigr) \overline{U}(x).
\end{equation}

Now we start to estimate the tail probability in (\ref{bound}). Choosing
some large $n$ such that $\sum_{i=n+1}^{\infty}1/i^{2}\leq1$.
Clearly, for
all $x>x_{0}$,
%
%
\begin{eqnarray}
\label{two terms} \Pr \Biggl( \sum_{i=n+1}^{\infty}X_{i}
\prod_{j=1}^{i}Y_{j}>x \Biggr) &
\leq& \Pr \Biggl( \sum_{i=n+1}^{\infty
}X_{i}
\prod_{j=1}^{i}Y_{j}>\sum
_{i=n+1}^{\infty}\frac
{x}{i^{2}} \Biggr)
\nonumber
\\
& \leq&\sum_{i=n+1}^{\infty}\Pr \Biggl(
X_{i}\prod_{j=1}^{i}Y_{j}>
\frac
{x}{%
i^{2}} \Biggr)
\nonumber
\\[-8pt]
\\[-8pt]
& \leq& \biggl( \sum_{i>\sqrt{x/x_{0}}}+\sum
_{n<i\leq\sqrt
{x/x_{0}}} \biggr) \Pr \Biggl( X^{\ast}\prod
_{j=1}^{i}Y_{j}^{\ast}>
\frac
{x}{i^{2}} \Biggr)
\nonumber
\\
& =&I_{1}(x)+I_{2}(n,x),
\nonumber
\end{eqnarray}
where $I_{2}(n,x)$ is understood as $0$ in case $n+1>\sqrt{x/x_{0}}$. First
we deal with $I_{1}(x)$. By Chebyshev's inequality,
\[
I_{1}(x)\leq x^{-\alpha}\mathrm{E} \bigl( X^{\ast} \bigr)
^{\alpha
}\sum_{i>%
\sqrt{x/x_{0}}}i^{2\alpha} \bigl(
\mathrm{E} \bigl( Y^{\ast} \bigr) ^{\alpha
} \bigr) ^{i}.
\]
This means that $I_{1}(x)$ converges to $0$ at least semi-exponentially fast
since $\mathrm{E} ( Y^{\ast} ) ^{\alpha}<1$. Thus,
%
%
\begin{equation}
\label{inf26} \lim_{x\rightarrow\infty}\frac{I_{1}(x)}{\overline{U}(x)}=0.
\end{equation}
Next we deal with $I_{2}(n,x)$. We further decompose it into three
parts as
%
%
\begin{eqnarray}
\label{three terms} I_{2}(n,x)& =&\sum_{n<i\leq\sqrt{x/x_{0}}}
\Pr \Biggl( X^{\ast
}\prod_{j=1}^{i}Y_{j}^{\ast}>
\frac{x}{i^{2}},0<X^{\ast}\leq\frac
{x}{i^{2}}%
 \Biggr)
\nonumber
\\
&&{} +\sum_{n<i\leq\sqrt{x/x_{0}}}\Pr \Biggl( X^{\ast}>
\frac
{x}{i^{2}}%
,\prod_{j=1}^{i}Y_{j}^{\ast}>1
\Biggr)
\nonumber
\\[-8pt]
\\[-8pt]
&&{} +\sum_{n<i\leq\sqrt{x/x_{0}}}\Pr \Biggl( X^{\ast
}\prod
_{j=1}^{i}Y_{j}^{\ast}>
\frac{x}{i^{2}},\prod_{j=1}^{i}Y_{j}^{\ast
}
\leq1 \Biggr)
\nonumber
\\
& =&I_{21}(n,x)+I_{22}(n,x)+I_{23}(n,x).
\nonumber
\end{eqnarray}
By conditioning on $X^{\ast}$ and then applying (\ref{inf1})--(\ref
{inf22}), we obtain
\begin{eqnarray*}
I_{21}(n,x)& \leq& K_{2}\sum_{n<i\leq\sqrt{x/x_{0}}}
\bigl( \mathrm {E} \bigl( Y^{\ast} \bigr) ^{\alpha}+\varepsilon
\bigr) ^{i}\Pr \biggl( X^{\ast
}Y^{\ast}>
\frac{x}{i^{2}} \biggr)
\\
& \sim& K_{1}K_{2} \bigl( \mathrm{E} \bigl(
X^{\ast} \bigr) ^{\alpha
}+\mathrm{%
E} \bigl(
Y^{\ast} \bigr) ^{\alpha} \bigr) \sum_{n<i\leq\sqrt{x/x_{0}}
}
\bigl( \mathrm{E} \bigl( Y^{\ast} \bigr) ^{\alpha}+\varepsilon
\bigr) ^{i}%
\overline{U} \biggl( \frac{x}{i^{2}} \biggr)
\\
& \leq&(1+\varepsilon)K_{1}K_{2} \bigl( \mathrm{E} \bigl(
X^{\ast
} \bigr) ^{\alpha}+\mathrm{E} \bigl( Y^{\ast} \bigr)
^{\alpha} \bigr) \overline {U}%
 ( x ) \sum
_{n<i\leq\sqrt{x/x_{0}}}i^{2(\alpha
+\varepsilon
)} \bigl( \mathrm{E} \bigl(
Y^{\ast} \bigr) ^{\alpha}+\varepsilon \bigr) ^{i}.
\end{eqnarray*}
Since $\mathrm{E} ( Y^{\ast} ) ^{\alpha}+\varepsilon<1$, it
follows that
%
%
\begin{equation}
\label{inf23} \lim_{n\rightarrow\infty}\limsup_{x\rightarrow\infty}
\frac
{I_{21}(n,x)}{%
\overline{U}(x)}=0.
\end{equation}
Applying both (\ref{inf1}) and (\ref{inf22}), we have
\[
I_{22}(n,x)\leq(1+\varepsilon)K_{1}K_{2}
\overline{G^{\ast}} ( 1 ) \overline{U} ( x ) \sum
_{n<i\leq\sqrt{x/x_{0}}%
}i^{2(\alpha+\varepsilon)} \bigl( \mathrm{E} \bigl(
Y^{\ast} \bigr) ^{\alpha}+\varepsilon \bigr) ^{i},
\]
which implies that
%
%
\begin{equation}
\label{inf24} \lim_{n\rightarrow\infty}\limsup_{x\rightarrow\infty}
\frac
{I_{22}(n,x)}{%
\overline{U}(x)}=0.
\end{equation}
Similarly, applying (\ref{inf22}) twice,
\begin{eqnarray*}
I_{23}(n,x)& \leq& K_{1}\sum_{n<i\leq\sqrt{x/x_{0}}}
\int_{0}^{1}\overline{U}%
 \biggl(
\frac{x}{i^{2}y} \biggr) \Pr \Biggl( \prod_{j=1}^{i}Y_{j}^{\ast
}
\in \mathrm{d}y \Biggr)
\\
& \leq&(1+\varepsilon)K_{1}\sum_{n<i\leq\sqrt{x/x_{0}}}
\overline {U} \biggl( \frac{x}{i^{2}} \biggr) \bigl( \mathrm{E} \bigl(
Y^{\ast} \bigr) ^{\alpha
-\varepsilon} \bigr) ^{i}
\\
& \leq&(1+\varepsilon)^{2}K_{1}\overline{U} ( x ) \sum
_{n<i\leq
\sqrt{x/x_{0}}}i^{2(\alpha+\varepsilon)} \bigl( \mathrm{E} \bigl(
Y^{\ast
} \bigr) ^{\alpha-\varepsilon} \bigr) ^{i},
\end{eqnarray*}
which, together with $\mathrm{E} ( Y^{\ast} ) ^{\alpha
-\varepsilon}<1$, gives that
%
%
\begin{equation}
\label{inf25} \lim_{n\rightarrow\infty}\limsup_{x\rightarrow\infty}
\frac
{I_{23}(n,x)}{%
\overline{U}(x)}=0.
\end{equation}
A combination of relations (\ref{two terms})--(\ref{inf25}) completes the
proof.
\end{pf}

\section{Proofs}\label{sec4}

\subsection{Proof of Theorem \texorpdfstring{\protect\ref{No1}(a)}{2.1(a)}}
\label{proof of (a)}

We first prove relation (\ref{r1}). It is easy to verify that
%
%
\begin{equation}
\label{M1} M_{n}\overset{\mathrm{d}} {=} ( X_{n}+M_{n-1}
) _{+}Y_{n},\qquad n\in\mathbb{N},
\end{equation}
where $\overset{\mathrm{d}}{=}$ denotes equality in distribution; see also
Theorem~2.1 of Tang and Tsitsiashvili \cite{TanTsi03}. We
proceed with
induction. For
$n=1$, it follows from Lemma~\ref{C2} that
%
%
\begin{equation}
\label{cited} \Pr ( M_{1}>x ) =\Pr ( X_{1,+}Y_{1}>x
) \sim \mu _{\alpha}\overline{F}(x)+\mathrm{E}X_{+}^{\alpha}
\overline{G}(x)=A_{1} \overline{F}(x)+B_{1}\overline{G}(x).
\end{equation}
Thus, relation (\ref{r1}) holds for $n=1$. Now we assume by induction that
relation (\ref{r1}) holds for $n-1\geq1$ and prove it for $n$. By this
induction assumption and Assumption~\ref{Assumption A}, we know that every
convex combination of the distributions of $X_{n}$ and $M_{n-1}$
belongs to
the class $\mathcal{R}_{-\alpha}^{\ast}\subset\mathcal{S}(0)$. Applying
Lemma~\ref{L2} with $\alpha=0$, we have
\[
\Pr ( X_{n}+M_{n-1}>x ) \sim ( 1+A_{n-1} ) \overline
{F}%
(x)+B_{n-1}\overline{G}(x),
\]
which, together with Assumption~\ref{Assumption A}, implies that every
convex combination of the distributions of $X_{n}+M_{n-1}$ and $Y_{n}$
belongs to the class $\mathcal{R}_{-\alpha}^{\ast}$. Applying Lemma~\ref{C2}, we obtain
\begin{eqnarray*}
\Pr ( M_{n}>x ) & =&\Pr \bigl( ( X_{n}+M_{n-1} )
_{+}Y_{n}>x \bigr)
\\
& \sim&\mu_{\alpha}\Pr ( X_{n}+M_{n-1}>x ) +\mathrm {E}
( X_{n}+M_{n-1} ) _{+}^{\alpha}
\overline{G}(x)
\\
& \sim& A_{n}\overline{F}(x)+B_{n}\overline{G}(x).
\end{eqnarray*}
Therefore, relation (\ref{r1}) holds for $n$.

Next we turn to relation (\ref{r2}). Introduce a sequence of random
variables $\{T_{n};n\in\mathbb{N}\}$ through the recursive equation
%
%
\begin{equation}
\label{T recursive equation} T_{n}= ( X_{n}+T_{n-1} )
Y_{n},\qquad n\in\mathbb{N},
\end{equation}
equipped with $T_{0}=0$. It is easy to see that $S_{n}\overset{\mathrm
{d}}{=}%
T_{n}$ for $n\in\mathbb{N}$. Then the proof of relation (\ref{r2})
can be
done by using the recursive equation (\ref{T recursive equation}) and going
along the same lines as in the proof of relation (\ref{r1}) above.

\subsection{Proof of Theorem \texorpdfstring{\protect\ref{No1}(b)}{2.1(b)}}
\label{proof of (b)}

Note that $A_{n}$ and $B_{n}$ increasingly converge to the finite
constants $%
A_{\infty}$ and $B_{\infty}$. Also recall Lemma~\ref{infinite}.
Hence, for
every $\delta>0$, there is some large integer $n_{0}$ such that both
%
%
\begin{equation}
\label{t3} ( A_{\infty}-A_{n_{0}} ) + ( B_{\infty
}-B_{n_{0}}
) \leq\delta
\end{equation}
and
%
%
\begin{equation}
\label{t4} \Pr \Biggl( \sum_{i=n_{0}+1}^{\infty}X_{i,+}
\prod_{j=1}^{i}Y_{j}>x \Biggr)
\lesssim\delta \bigl( \overline{F}(x)+\overline{G}(x) \bigr)
\end{equation}
hold. Now we start to deal with $\Pr ( M_{\infty}>x ) $.
On the
one hand, for every $\varepsilon>0$, by Theorem~\ref{No1}(a), relation
(\ref{t4}), and Assumption~\ref{Assumption A}, in turn, we obtain
%
%
\begin{eqnarray}
\label{t5} \Pr ( M_{\infty}>x ) & \leq&\Pr \bigl( M_{n_{0}}> ( 1-
\varepsilon ) x \bigr) +\Pr \Biggl( \sum_{i=n_{0}+1}^{\infty
}X_{i,+}
\prod_{j=1}^{i}Y_{j}>\varepsilon x
\Biggr)
\nonumber
\\
& \lesssim& A_{n_{0}}\overline{F}\bigl( ( 1-\varepsilon ) x
\bigr)+B_{n_{0}}%
\overline{G}\bigl( ( 1-\varepsilon ) x\bigr)+
\delta \bigl( \overline{F}%
(\varepsilon x)+\overline{G}(\varepsilon x)
\bigr)
\nonumber
\\[-8pt]
\\[-8pt]
& \sim& ( 1-\varepsilon ) ^{-\alpha} \bigl( A_{n_{0}}\overline
{F}%
(x)+B_{n_{0}}\overline{G}(x) \bigr) +\delta
\varepsilon^{-\alpha
} \bigl( \overline{F}(x)+\overline{G}(x) \bigr)
\nonumber
\\
& \leq& \bigl( ( 1-\varepsilon ) ^{-\alpha}A_{\infty
}+\delta
\varepsilon^{-\alpha} \bigr) \overline{F}(x)+ \bigl( ( 1-\varepsilon )
^{-\alpha}B_{\infty}+\delta\varepsilon^{-\alpha} \bigr)
\overline{G}(x).
\nonumber
\end{eqnarray}
On the other hand, by Theorem~\ref{No1}(a) and relation (\ref{t3}),
%
%
\begin{equation}
\label{t6} \Pr ( M_{\infty}>x ) \geq\Pr ( M_{n_{0}}>x ) \gtrsim (
A_{\infty}-\delta ) \overline{F}(x)+ ( B_{\infty
}-\delta )
\overline{G}(x).
\end{equation}
By the arbitrariness of $\delta$ and $\varepsilon$ in (\ref{t5}) and
(\ref{t6}), we obtain relation (\ref{r1}) for $n=\infty$.

\subsection{Proof of Theorem \texorpdfstring{\protect\ref{No1}(c)}{2.1(c)}}
\label{proof of (c)}

First we establish an asymptotic upper bound for $\Pr (
S_{\infty
}>x ) $. As in the proof of Theorem~\ref{No1}(b), for every $%
\delta>0$, suitably choose some large integer $n_{0}$ such that
relations (\ref{t3}), (\ref{t4}), and the relation
%
%
\begin{equation}
\label{t7} -\delta\leq C_{\infty}-C_{n_{0}}\leq\delta
\end{equation}
hold simultaneously. For every $\varepsilon>0$, by Theorem~\ref{No1}(a),
relation (\ref{t4}), Assumption~\ref{Assumption A}, and relation
(\ref{t7}),
in turn, we obtain
\begin{eqnarray*}
\Pr ( S_{\infty}>x ) & \leq&\Pr \bigl( S_{n_{0}}> ( 1-\varepsilon )
x \bigr) +\Pr \Biggl( \sum_{i=n_{0}+1}^{\infty
}X_{i,+}
\prod_{j=1}^{i}Y_{j}>\varepsilon x
\Biggr)
\\
& \lesssim& \bigl( A_{n_{0}}\overline{F}\bigl( ( 1-\varepsilon ) x
\bigr)+C_{n_{0}}\overline{G}\bigl( ( 1-\varepsilon ) x\bigr) \bigr) +\delta
\bigl( \overline{F}(\varepsilon x)+\overline{G}(\varepsilon x) \bigr)
\\
& \sim& ( 1-\varepsilon ) ^{-\alpha} \bigl( A_{n_{0}}\overline
{F}%
(x)+C_{n_{0}}\overline{G}(x) \bigr) +\delta
\varepsilon^{-\alpha
} \bigl( \overline{F}(x)+\overline{G}(x) \bigr)
\\
& \leq& \bigl( ( 1-\varepsilon ) ^{-\alpha}A_{\infty
}+\delta
\varepsilon^{-\alpha} \bigr) \overline{F}(x)+ \bigl( ( 1-\varepsilon )
^{-\alpha} ( C_{\infty}+\delta ) +\delta \varepsilon ^{-\alpha}
\bigr) \overline{G}(x).
\end{eqnarray*}
Since $\delta$ and $\varepsilon$ are arbitrary positive constants, it
follows that
\[
\Pr ( S_{\infty}>x ) \lesssim A_{\infty}\overline{F}%
(x)+C_{\infty}
\overline{G}(x).
\]

For the corresponding asymptotic lower bound, as analyzed in the proof of
Theorem~\ref{No1}(a), it suffices to prove that
%
%
\begin{equation}
\label{lower} \Pr ( T_{\infty}>x ) \gtrsim A_{\infty}\overline
{F}(x)+C_{\infty
}\overline{G}(x),
\end{equation}
where $T_{\infty}$ is the weak limit of the sequence $\{T_{n};n\in
\mathbb{N%
}\}$ defined by (\ref{T recursive equation}). We apply the method developed
by Grey \cite{Gre94} to prove (\ref{lower}).
Consider the stochastic difference
equation
%
%
\begin{equation}
\label{stochastic difference equation} T_{\infty}\overset{\mathrm{d}} {=} ( X+T_{\infty} ) Y,
\end{equation}
which inherits a stochastic structure from (\ref{T recursive equation}).
Note that the weak solution of (\ref{stochastic difference equation}) exists
and is unique. Furthermore, the limit distribution of $T_{n}$ is identical
to this unique solution and, hence, it does not depend on the starting
random variable $T_{0}$. See Vervaat \cite{Ver79}
and Goldie \cite{Gol91} for these and
related statements.

It is easy to check that $q=\Pr ( T_{\infty}>0 ) >0$; see the
proof of Theorem~1 of Grey \cite{Gre94} for a
similar argument. Construct a new
starting random variable $\tilde{T}_{0}$ independent of $\{X_{1},X_{2},
\ldots;Y_{1},Y_{2},\ldots\}$ with tail
%
%
\begin{equation}
\label{starting random variable} \Pr ( \tilde{T}_{0}>x ) =q\Pr ( XY>x ) \mathbf
{1}%
_{(x\geq0)}+\Pr ( T_{\infty}>x ) \mathbf{1}_{(x<0)}.
\end{equation}
Starting with $\tilde{T}_{0}$, the recursive equation (\ref{T recursive
equation}) generates the sequence $\{\tilde{T}_{n};n\in\mathbb{N}\}$
correspondingly. Comparing (\ref{starting random variable}) with
(\ref{stochastic difference equation}), we see that $\tilde{T}_{0}$ and, hence,
every $\tilde{T}_{n}$ are stochastically not greater than $T_{\infty}$;
namely, it holds for all $x\in\mathbb{R}$ and all $n\in\{0\}\cup
\mathbb{N}
$ that
%
%
\begin{equation}
\label{lower 1} \Pr ( T_{\infty}>x ) \geq\Pr ( \tilde {T}_{n}>x ).
\end{equation}
Furthermore, it holds that
\[
\Pr ( \tilde{T}_{0}>x ) \sim q\Pr ( X_{+}Y>x ) \sim q\mu
_{\alpha}\overline{F}(x)+q\mathrm{E}X_{+}^{\alpha}
\overline{G}(x),
\]
where the last step is analogous to (\ref{cited}). Thus, by Assumption~\ref{Assumption A}, the distribution of $\tilde{T}_{0}$ belongs to the
class $%
\mathcal{R}_{-\alpha}^{\ast}$. Then, by going along the same lines
of the
proof of Theorem~\ref{No1}(a) and using equation (\ref{T recursive
equation}) starting with $\tilde{T}_{0}$, we obtain
%
%
\begin{equation}
\label{lower 2} \Pr ( \tilde{T}_{n}>x ) \sim\tilde{A}_{n}
\overline {F}(x)+\tilde{C}%
_{n}\overline{G}(x)
\end{equation}
with
\[
\tilde{A}_{n}=\sum_{i=1}^{n}
\mu_{\alpha}^{i}+q\mu_{\alpha
}^{n+1},\qquad
\tilde{C}_{n}=\sum_{i=1}^{n}
\mu_{\alpha}^{i-2}\mathrm{E}\tilde{T}%
_{n-i+1,+}^{\alpha}+q
\mu_{\alpha}^{n}\mathrm{E}X_{+}^{\alpha}.
\]
Since $\tilde{T}_{n}$ weakly converges to $T_{\infty}\overset
{\mathrm
{d}}{=}%
S_{\infty}$ and $\mu_{\alpha}<1$, it is easy to see that $%
\lim_{n\rightarrow\infty}\tilde{A}_{n}=A_{\infty}$ and $%
\lim_{n\rightarrow\infty}\tilde{C}_{n}=C_{\infty}$, with the latter
subject to a straightforward application of the dominated convergence
theorem. Thus, substituting (\ref{lower 2}) into (\ref{lower 1}) and letting
$n\rightarrow\infty$ on the right-hand side of the resulting formula, we
arrive at relation (\ref{lower}) as desired.

\subsection{Sketch of the proof of Corollary \texorpdfstring{\protect\ref{No2}}{2.1}}
\label{proof of cor}

Clearly, the recursive equations (\ref{M1}), (\ref{T recursive equation}),
and the identity $S_{n}\overset{\mathrm{d}}{=}T_{n}$ for $n\in
\mathbb{N}$
still hold since $\{(X_{1},Y_{1}),(X_{2},Y_{2}),\ldots\}$ is a
sequence of
i.i.d. random pairs. Introduce four independent random variables
$X^{\prime
} $, $\check{X}^{\prime}$, $Y^{\prime}$, and $\check{Y}^{\prime}$ with
distributions $F$, $F^{2}$, $G$, and $G^{2}$, respectively, and let
them be
independent of $\{(X_{1},Y_{1}),(X_{2},Y_{2}),\ldots\}$. Using
decomposition (\ref{Pi}), we have
%
%
\begin{eqnarray}
\label{decomposition} \Pr ( M_{n}>x ) & =&\Pr \bigl( ( X_{n}+M_{n-1}
) _{+}Y_{n}>x \bigr)
\nonumber
\\
& =&(1+\theta)\Pr \bigl( \bigl( X^{\prime}+M_{n-1} \bigr)
_{+}Y^{\prime
}>x \bigr) -\theta\Pr \bigl( \bigl(
\check{X}^{\prime
}+M_{n-1} \bigr) _{+}Y^{\prime}>x
\bigr)
\\
&&{} -\theta\Pr \bigl( \bigl( X^{\prime}+M_{n-1} \bigr)
_{+}\check {Y}%
^{\prime}>x \bigr) +\theta\Pr \bigl(
\bigl( \check{X}^{\prime
}+M_{n-1} \bigr) _{+}
\check{Y}^{\prime}>x \bigr).
\nonumber
\end{eqnarray}
When $n=1$, applying Lemma~\ref{C2} to each term on the right-hand side
of (\ref{decomposition}) gives
%
%
\begin{equation}
\label{c1} \Pr ( M_{1}>x ) =\Pr ( X_{1,+}Y_{1}>x
) \sim A_{1}^{\prime}\overline{F}(x)+B_{1}^{\prime}
\overline{G}(x).
\end{equation}
Then, as in the proof of Theorem~\ref{No1}(a), proceeding with induction
according to (\ref{decomposition}) leads to (\ref{r3}). Relation
(\ref{r4})
can be derived similarly. This proves Corollary~\ref{No2}(a).

Corollary~\ref{No2}(b), (c) can be verified by the similar ideas used in
proving Theorem~\ref{No1}(b),~(c). The key ingredient is establishing a
relation similar to (\ref{bound}). Write $Z=XY$, $Z_{1}=X_{1}Y_{1}$, $%
Z_{2}=X_{2}Y_{2}$, and so on. It follows from (\ref{c1}) that
\[
\Pr ( Z>x ) +\overline{G}(x)\asymp\overline {F}(x)+\overline{G}(x).
\]
As in the proof of Lemma~\ref{infinite}, we can construct independent random
variables $Z^{\ast}$ and $Y^{\ast}$ both with tails equal to
$K_{1} (
\overline{F}(x)+\overline{G}(x) ) $ for all large $x$ such that
$Z\leq
_{\mathrm{st}}Z^{\ast}$, $Y\leq_{\mathrm{st}}Y^{\ast}$, and
$\mathrm
{E}%
 ( Y^{\ast} ) ^{\alpha}<1$. For some large $n$ such that $%
\sum_{i=n+1}^{\infty}1/i^{2}\leq1$, we write
\begin{eqnarray*}
\Pr \Biggl( \sum_{i=n+1}^{\infty}X_{i}
\prod_{j=1}^{i}Y_{j}>x \Biggr) &
\leq& \Pr \Biggl( \sum_{i=n+1}^{\infty
}X_{i}
\prod_{j=1}^{i}Y_{j}>\sum
_{i=n+1}^{\infty}\frac
{x}{i^{2}} \Biggr)
\\
& =&\sum_{i=n+1}^{\infty}\Pr \Biggl(
Z_{i}\prod_{j=1}^{i-1}Y_{j}>
\frac
{x}{%
i^{2}} \Biggr)
\\
& \leq&\sum_{i=n+1}^{\infty}\Pr \Biggl(
Z^{\ast
}\prod_{j=1}^{i-1}Y_{j}^{\ast}>
\frac{x}{i^{2}} \Biggr).
\end{eqnarray*}
Then, going along the same lines of the rest of the proof of Lemma~\ref{infinite}, we obtain
\[
\lim_{n\rightarrow\infty}\limsup_{x\rightarrow\infty}\frac
{1}{\overline{F%
}(x)+\overline{G}(x)}\Pr
\Biggl( \sum_{i=n+1}^{\infty
}X_{i}\prod
_{j=1}^{i}Y_{j}>x \Biggr) =0,
\]
which suffices for our purpose.


\section*{Acknowledgements}

The authors are very grateful to the two
reviewers and the associate editor for their insightful comments and
constructive suggestions, which have helped significantly improve this work.
Li's research was supported by the National Natural Science Foundation of
China (Grant No.~11201245) and the Research Fund for the Doctoral
Program of
Higher Education of China (Grant No.~20110031120003). Tang's research was
supported by a Centers of Actuarial Excellence (CAE) Research Grant
(2013--2016) from the Society of Actuaries.


%

\printhistory
\end{document}